\documentclass{amsart}
\usepackage[latin1]{inputenc}
\usepackage[T1]{fontenc}
\usepackage{amsfonts}
\usepackage{amsmath}
\usepackage{amssymb}
\usepackage{pstricks}
\usepackage{pst-grad}
\usepackage{multido}
\usepackage{amsthm}

\newcommand{\cA}{\mathcal{A}}

\newcommand{\cF}{\mathcal{F}}

\newcommand{\cH}{\mathcal{H}}

\newcommand{\cM}{\mathcal{M}}

\newcommand{\sg}{\mathcal{h}}
\newcommand{\sd}{\mathcal{i}}

\newcommand{\nZ}{\mathbb{Z}}

\newcommand{\al}{\alpha}

\newcommand{\si}{\sigma}
\newcommand{\la}{\lambda}

\newcommand{\va}{\varphi}

\newcommand{\tA}{\tilde{A}}

\newcommand{\tC}{\tilde{C}}

\newcommand{\ov}{\overline}

\newcommand{\Hi}{H^{(i)}}
\newcommand{\Hj}{H^{(j)}}

\newtheorem{Th}{Theorem}[section]
\newtheorem{Lem}[Th]{Lemma}
\newtheorem{Cor}[Th]{Corollary}
\newtheorem{Prop}[Th]{Proposition}
\newtheorem{Def-Prop}[Th]{Definition-Proposition}

\theoremstyle{definition}

\newtheorem{what}[Th]{}

\theoremstyle{remark}

\begin{document}

\title{On the lowest two-sided cell in affine Weyl groups}
\author{J\'er\'emie Guilhot}
\address{School of Mathematics and Statistics F07\\
University of Sydney NSW 2006\\
Australia} 
\email{guilhot@maths.usyd.edu.au}

\begin{abstract}
Bremke and Xi determined the lowest two-sided cell for affine Weyl groups with unequal parameters and showed that it consists of at most $|W_{0}|$ left cells where $W_{0}$ is the associated finite Weyl group. We prove that this bound is exact. Previously, this was known in the equal parameter case and when the parameters were coming from a graph automorphism. Our argument uniformly works for any choice of parameters.
\end{abstract}

\maketitle


\section{Introduction}
This paper is concerned with the theory of Kazhdan-Lusztig cells in a Coxeter group $W$, following the general setting of Lusztig \cite{bible}. This involves a weight function $L$ which is an integer-valued function on $W$ such that $L(ww')=L(w)+L(w')$ whenever $\ell(ww')=\ell(w)+\ell(w')$  ($\ell$ is the usual length function on $W$). We shall only consider weight function such that $L(w)>0$ for all $w\neq 1$. The case where $L=\ell$ is known as the equal parameter case.

The partition of $W$ into cells is known to play an important role in the study of the representations of the corresponding Hecke algebra.

Here, we are primarily concerned with affine Weyl groups. Let $W$ be an irreducible affine Weyl group, together with a weight function $L$. Let $W_{0}$ be the finite Weyl group associated to $W$. 

In the equal parameter case, Shi \cite{Shi1} described the lowest two-sided cell $c_{0}$ with respect to the preorder $\leq_{LR}$ using the Lusztig $a$-function (see \cite{bible} for further details on the $a$-function) as follows
$$c_{0}=\{w\in W|\ a(w)=\ell(w_{0})\}$$
where $w_{0}$ is the longest element of $W_{0}$. Then, he gave an upper bound for the number of left cells contained in $c_{0}$, namely $|W_{0}|$. In \cite{Shi2}, he showed that this bound is exact and he described the left cells in $c_{0}$. His proof involved some deep properties of the Kazhdan-Lusztig polynomials in the equal parameter case, such as the positivity of the coefficients.

Using the positivity property, Lusztig ([6]) proved, in the equal parameter case,  a number of results concerning the $a$-function, such as
\begin{equation*}
\text{if $z\leq_{LR}z'$ then $a(z)\geq a(z')$}  \tag{P4}
\end{equation*}
and
\begin{equation*}
 \text{if } z\leq_{LR} z' \text{ and } a(z)=a(z') \text{ then } z \sim_{LR} z'.  \tag{P11}
\end{equation*}
Since we know that $a(z)\leq \ell(w_{0})$ for any $z\in W$, this shows that $c_{0}$ is the lowest two-sided cell with respect to $\leq_{LR}$.

In the unequal parameter case, the positivity of the coefficients of the Kazhdan-Lusztig polynomials does not hold anymore. However, Bremke \cite{Bremke} and Xi \cite{Xi} proved that the lowest two-sided cell $c_{0}$ can be described in the same way as in \cite{Shi1}, using the general $a$-function. Let $J$ be a subset of $S$ and consider the corresponding parabolic subgroup $W_{J}=\sg J\sd$. We denote by $w_{J}$ the longest element in $W_{J}$. Let
$$\displaystyle \tilde{\nu}=\max_{J\subset S, W_{J}\simeq W_{0}}L(w_{J}).$$
Then, we have
$$c_{0}=\{w\in W|\ a(w)=\tilde{\nu}\}.$$
Let $\textbf{S}$ be the set which consists of all the subset $J$ of $S$ such that $W_{J}\simeq W_{0}$ and $L(w_{J})=\tilde{\nu}$. Then we have the following alternative
description of $c_0$, and it is this description with which we work throughout the paper
$$c_{0}=\{w\in W|\ w=x.w_{J}.y,\ x,y\in W,\ J\in \textbf{S}\}$$
where for any $x,y,w\in W$ the notation $w=x.y$ means that $w=xy$ and $\ell(w)=\ell(x)+\ell(y)$ (and similarly $w=x.y.z$ for $w,x,y,z\in W$).

In \cite{Bremke}, Bremke showed that $c_{0}$ contains at most $|W_{0}|$ left cells. She proved that this bound is exact when the parameters are coming from a graph automorphism and, again, the proof involved some deep properties of the Kazhdan-Lusztig polynomials in the equal parameter case.

In this paper, we will prove that $|W_{0}|$ is the exact bound for the number of left cells in $c_{0}$ for any choice of parameters; see Theorem \ref{ltsc}. The method is based on a variation of the induction of left cells as in \cite{Geck}. The main new ingredient is a geometric argument to find a ``local'' bound on the degree of the structure constants of the Hecke algebra of $W$ with respect to the standard basis; see Theorem \ref{bound}. Our proof uniformly works for any choice of parameters.

\section{Multiplication of the standard basis and geometric realization}
In this section, we introduce the Hecke algebra of a Coxeter group with respect to a weight function. Then, we present a geometric realization of an affine Weyl group. Finally, using this geometric realization, we give a bound on the degree of the structure constants with respect to the standard basis.
\begin{what}{\bf Weight functions, Hecke algebras.}
In this section, $(W,S)$ denotes an arbitrary Coxeter system. The basic reference is \cite{bible}. Let $L$ be a weight function. In this paper, we will only consider the case where $L(w)>0$ for all $w\neq 1$. A weight function is completely determined by its value on $S$ and must only satisfy $L(s)=L(t)$ if $s$ and $t$ are conjugate.

Let $\cA=\nZ[v,v^{-1}]$ and $\cH$ be the generic Iwahori-Hecke algebra associated to $(W,S)$ with parameters $\{L(s)\mid s\in S\}$. $\cH$ has an $\cA$-basis $\{T_{w}\mid w\in W\}$, called the standard basis, with multiplication given by
\begin{equation*}
T_{s}T_{w}=
\begin{cases}
T_{sw}, & \mbox{if } sw>w,\\
T_{sw}+(v^{L(s)}-v^{-L(s)})T_{w}, &\mbox{if } sw<w,
\end{cases}
\end{equation*}
(here, ``<'' denotes the Bruhat order) where $s\in S$ and $w\in W$.

Let $x,y\in W$. We write
$$\displaystyle T_{x}T_{y}=\sum_{z\in W}f_{x,y,z}T_{z}$$
where $f_{x,y,z}\in\cA$ are the structure constants with respect to the standard basis.

In this paper, we will be mainly interested in the case where $W$ is an irreducible affine Weyl group (with corresponding Weyl group $W_{0}$). In that case, it is known that there is a global bound for the degrees of the structure constants $f_{x,y,z}$. Namely, set $\nu=\ell(w_{0})$, $\tilde{\nu}=L(w_{0})$, where $w_{0}$ is the longest element of $W_{0}$, and $\xi_{s}=v^{L(s)}-v^{-L(s)}$ for $s\in S$. In \cite{Bremke}, Bremke proved that
\begin{enumerate}
\item As a polynomial in $\xi_{s}$, $s\in S$, the degree of $f_{x,y,z}$ is at most $\nu$
\item The degree of $f_{x,y,z}$ in $v$ is at most $\tilde{\nu}$.
\end{enumerate}

Our aim will be to find a ``local'' bound for the degrees of these polynomials, which depends on $x,y \in W$. For this purpose, we will work with a geometric realization of $W$, as described in the next section. 
\end{what}

\begin{what}{\bf Geometric realization.}
\label{gem}
In this section, we present a geometric realization of an affine Weyl group. The basic references are \cite{Bremke,Lus1,Xi}. \\

Let $V$ be an euclidean space of finite dimension $r\geq 1$. Let $\Phi$ be an irreducible root system of rank $r$ and $\check{\Phi}\subset V^{*}$ the dual root system. We denote the coroot corresponding to $\alpha\in\Phi$ by $\check{\alpha}$ and we write $\mathcal{h}x,y\mathcal{i}$ for the value of $y\in V^{*}$ at $x\in V$. Fix a set of positive roots $\Phi^{+}\subset \Phi$. Let $W_{0}$ be the Weyl group of $\Phi$. For $\alpha\in\Phi^{+}$ and $n\in \mathbb{Z}$, we define a hyperplane
$$H_{\alpha,n}=\{x\in V\mid \mathcal{h}x,\check{\alpha}\mathcal{i}=n\}.$$
Let
$$\mathcal{F}=\{H_{\alpha,n}\mid \alpha\in \Phi^{+}, n\in\mathbb{Z}\}.$$
Any $H\in\mathcal{F}$ defines an orthogonal reflection $\sigma_{H}$ with fixed point set $H$. We denote by $\Omega$ the group generated by all these reflections, and we regard $\Omega$ as acting on the right on $V$. An alcove is a connected component of the set
$$V-\underset{H\in\mathcal{F}}{\bigcup}H.$$
$\Omega$ acts simply transitively on the set of alcoves $X$.

Let $S$ be the set of $\Omega$-orbits in the set of faces (codimension 1 facets) of alcoves. Then $S$ consists of $r+1$ elements which can be represented as the $r+1$ faces of an alcove. If a face $f$ is contained in the orbit $t\in S$, we say that $f$ is of type $t$.

Let $s\in S$. We define an involution $A\rightarrow sA$ of $X$ as follows.  Let $A\in X$; then $sA$ is the unique alcove distinct from $A$ which shares with $A$ a face of type $s$.  The set of such maps generates a group of permutations of $X$ which is a Coxeter group $(W,S)$. In our case, it is the affine Weyl group usually denoted $\tilde{W_{0}}$. We regard $W$ as acting on the left on $X$. It acts simply transitively and commutes with the action of $\Omega$.

Let $L$ be a weight function on $W$. In \cite{Bremke}, Bremke showed that if a hyperplane $H$ in $\cF$ supports faces of type $s,t\in S$ then $s$ and $t$ are conjugate in $W$ which implies that $L(s)=L(t)$. Thus we can associate a weight $c_{H}\in\mathbb{Z}$ to $H\in \cF$ such that $c_{H}=L(s)$ if $H$ supports a face of type $s$. 

Assume that $W$ is not of type $\tilde{A}_{1}$ or $\tilde{C}_{r}$ ($r\geq 2$); Let $H,H'$ be two parallel hyperplanes such that $H$ supports a face of type $s$ and $H'$ a face of type $s'$. Then, Bremke \cite{Bremke} proved that $s$ and $s'$ are conjugate. In other words, if $W$ is not of type $\tilde{A}_{1}$ or $\tilde{C}_{r}$ then any two parallel hyperplanes have the same weight.

In this paper, we will often have to distinguish the case where $W$ is of type $\tilde{A}_{1}$ or $\tilde{C}_{r}$, because of this property.

 In the case where $W$ is of type $\tilde{C}_{r}$ with generators $s_{1},...,s_{r+1}$ and $W_{0}$ is generated by $s_{1},...,s_{r}$, by symmetry of the Dynkin diagram, we can assume that $c_{s_{1}}\geq c_{s_{r+1}}$. Similarly, if $W$ is of type $\tilde{A}_{1}$ with generators $s_{1},s_{2}$ and $W_{0}=\mathcal{h}s_{1}\mathcal{i}$, we can assume that $c_{s_{1}}\geq c_{s_{2}}$.

For a $0$-dimensional facet $\la$ of an alcove, define
$$m(\la)=\underset{H\in\cF,\ \la\in H}{\sum}c_{H}.$$
We say that $\la$ is a special point if $m(\la)$ is maximal.  

Let $T$ be the set of all special points. For $\la\in T$, denote by $W_{\la}$ the stabilizer of the set of alcoves containing $\la$ in their closure with respect to the action of $W$ on $X$. It is a maximal parabolic subgroup of $W$. Let $S_{\la}=S\cap W_{\la}$ and write $w_{\la}$ for the longest element of $W_{\la}$. Note that, following \cite{Bremke}, and with our convention for $\tilde{C_{r}}$ and $\tilde{A}_{1}$, $0\in V$ is a special point. Moreover, if $\la=0\in V$,  the definition of $W_{\lambda}$ is consistent with the definition of $W_{0}$ given before.

Let $\la$ be a special point, a quarter with vertex $\la$ is a connected component of
$$V-\underset{H,\ \la\in H}{\bigcup}H.$$
It is an open simplicial cone. It has $r$ walls.

Let $H=H_{\alpha,n}\in \cF$. Then $H$ divides $V-H$ into two half-spaces
\begin{align*}
V_{H}^{+}&=\{x\in V\mid \sg x,\check{\alpha}\sd>n\},\\
V_{H}^{-}&=\{x\in V\mid \sg x,\check{\alpha}\sd<n\}.
\end{align*}

Finally, let $A_{0}$ be the fundamental alcove defined by
$$A_{0}=\{x\in V\mid 0<\sg x,\check{\alpha}\sd<1 \text{ for all $\alpha\in\Phi^{+}$}\}.$$

Let $A\in X$, $w\in W$. It is well known that the length of $w$ is the number of hyperplanes which separate $A$ and $wA$. 
\end{what}

\begin{what}{\bf Multiplication of the standard basis.}
Let $(W,S)$ be an irreducible affine Weyl group associated to the Weyl group $(W_{0},S_{0})$. Recall that, for $x,y\in W$, we have
$$\displaystyle T_{x}T_{y}=\sum_{z\in W}f_{x,y,z}T_{z}.$$
After the preparations in 2.2, we will now be able to find a ``local'' bound for the degree of the polynomials $f_{x,y,z}$ which depends on $x,y\in W$.   

For two alcoves $A,B\in X$, let
$$H(A,B)=\{H\in\mathcal{F}\mid H \text{ separates } A \text{ and } B\}.$$
Let $\overline{\mathcal{F}}$ be the set of directions of hyperplanes in $\cF$. For $i\in \overline{\mathcal{F}}$, we denote by $\mathcal{F}_{i}$ the set of all hyperplanes $H\in\mathcal{F}$ of direction $i$. The connected components of
$$V-\underset{H\in\mathcal{F}_{i}}{\bigcup}H$$
are called ``strip of direction $i$''. We denote by $U_{i}(A)$ the unique strip of direction $i$ which contains $A$, for $A\in X$. There exists a unique $\alpha\in\Phi^{+}$ and a unique $n\in\mathbb{Z}$ such that
$$U_{i}(A)=\{\mu\in V\mid n<\mathcal{h}\mu,\check{\alpha}\mathcal{i}<n+1\},$$
in other words
$$U_{i}(A)=V_{H_{\alpha,n}}^{+}\cap V_{H_{\alpha,n+1}}^{-}.$$
We say that $U_{i}(A)$ is defined by $H_{\alpha,n}$ and $H_{\alpha,n+1}$.

Note that our definition of strips is slightly different from the one in \cite{Bremke}, where the strips were the connected components of 
$$V-\underset{ c_{H}=c_{i}}{\underset{H\in\mathcal{F}_{i}}{\bigcup}}H \quad\text{where}\quad c_{i}=\underset{H\in\mathcal{F},\ \overline{H}=i}{\max}c_{H}.
$$

In fact, as noticed in Section \ref{gem}, if $W$ is not of type $\tilde{A}_{1}$ or $\tilde{C}_{r}$, then the two definitions are the same.

Let $\sigma\in\Omega$. Let $i,j\in \overline{\mathcal{F}}$ such that $\sigma(i)=j$. We have
$$(U_{i}(A))\sigma=U_{j}(A\sigma)$$
and the strip $U_{j}(A\sigma)$ is defined by the two hyperplanes $(H_{\alpha,n_{i}})\sigma$ and $(H_{\alpha,n_{i}+1})\sigma$. \\
Let $x,y\in W$; then we define
\begin{align*}
H_{x,y}&=\{H\in \mathcal{F}\mid H\in H(A_{0},yA_{0})\cap H(yA_{0},xyA_{0})\},\\
I_{x,y}&=\{i\in\overline{\mathcal{F}}\mid\exists H, \overline{H}=i, H\in H_{x,y}\}.
\end{align*}
For $i\in I_{x,y}$, let
$$c_{x,y}(i)=\underset{H\in\mathcal{F},\ \overline{H}=i,\ H\in H_{x,y}}{\text{max }} c_{H}$$
and
$$c_{x,y}=\underset{i\in I_{x,y}}{\sum}c_{x,y}(i).$$

We are now ready to state the main result of this section.
\begin{Th}
\label{bound}
Let $x,y\in W$ and
$$T_{x}T_{y}=\underset{z\in W}{\sum}f_{x,y,z}T_{z}\ \ \text{ where $f_{x,y,z}\in\cA$}.$$
Then, the degree of $f_{x,y,z}$ in $v$ is at most $c_{x,y}$.
\end{Th}
\end{what}

\section{Proof of Theorem \ref{bound}}
In order to prove Theorem \ref{bound}, we will need the following lemmas.
\begin{Lem}
\label{Lem1}
Let $x,y\in W$ and $s\in S$ be such that $x<xs$ and $y<sy$. We have
$$c_{xs,y}=c_{x,sy}.$$
\end{Lem}
\begin{proof}
Let $H_{s}$ be the unique hyperplane which separates $yA_{0}$ and $syA_{0}$. Since $x<xs$ and $y<sy$, one can see that
\begin{align*}
H(A_{0},yA_{0})\cup\{H_{s}\}&=H(A_{0},syA_{0}),\\
H(A_{0},yA_{0})\cap\{H_{s}\}&=\emptyset,
\end{align*}
and
\begin{align*}
H(syA_{0},xsyA_{0})\cup\{H_{s}\}&=H(yA_{0},xsyA_{0}),\\
H(syA_{0},xsyA_{0})\cap\{H_{s}\}&=\emptyset.
\end{align*}
Therefore we have
\begin{align*}
H_{x,sy}&=H(A_{0},syA_{0})\cap H(syA_{0},xsyA_{0})\\
&= (H(A_{0},yA_{0})\cup\{H_{s}\})\cap H(syA_{0},xsyA_{0})\\
&= (H(A_{0},yA_{0})\cap H(syA_{0},xsyA_{0}))\cup (\{H_{s}\}\cap H(syA_{0},xsyA_{0}))\\
&=H(A_{0},yA_{0})\cap H(syA_{0},xsyA_{0})
\end{align*}
and
\begin{align*}
H_{xs,y}&=H(yA_{0},xsyA_{0})\cap H(A_{0},yA_{0})\\
&=(H(syA_{0},xsyA_{0})\cup\{H_{s}\})\cap H(A_{0},yA_{0})\\
&=(H(syA_{0},xsyA_{0})\cap H(A_{0},yA_{0}))\cup(\{H_{s}\}\cap H(A_{0},yA_{0}))\\
&=H(syA_{0},xsyA_{0})\cap H(A_{0},yA_{0})\\
&=H_{x,sy}.
\end{align*}
Thus $c_{x,sy}=c_{xs,y}$.
\end{proof}
\begin{Lem}
\label{lem0}
Let $x,y\in W$ and $s\in S$ be such that $xs<x$ and $sy<y$. We have
$$c_{xs,sy}\leq c_{x,y}.$$
\end{Lem}
\begin{proof}
Let $H_{s}$ be the unique hyperplane which separate $yA_{0}$ and $syA_{0}$. One can see that
$$H_{xs,sy}=H_{x,y}-\{H_{s}\}.$$
The result follows.
\end{proof}

\begin{Lem}
\label{Lem2}
Let $x,y\in W$ and $s\in S$ be such that $xs<x$ and $sy<y$. Let $H_{s}$ be the unique hyperplane which separates $yA_{0}$ and $syA_{0}$. Then we have
$$\overline{H_{s}}\notin I_{xs,y} \quad \text{and}\quad \overline{H_{s}}\in I_{x,y}.$$
\end{Lem}
\begin{proof}
We have
\begin{align*}
sy<y&\Longrightarrow H_{s}\in H(A_{0},yA_{0}),\\
xs<x&\Longrightarrow H_{s}\in H(yA_{0},xyA_{0}).
\end{align*}
Thus $H_{s}\in H_{x,y}$ and $\overline{H}_{s}\in I_{x,y}$.

Let $\alpha_{s}\in\Phi^{+}$ and $n_{s}\in\mathbb{Z}$ be such that $H_{s}=H_{\alpha_{s},n_{s}}$.
Assume that $n_{s}\geq 1$ (the case where $n_{s}\leq 0$ is similar).
Since $H_{s}\in H(A_{0},yA_{0})$ and $yA_{0}$ has a facet contained in $H_{s}$, we have
$$n_{s}<\mathcal{h}\mu,\check{\al}_{s}\mathcal{i}<n_{s}+1\text{ for all $\mu\in yA_{0}$}.$$
Therefore, for all $m>n_{s}$, we have $H_{\alpha_{s},m}\notin H(A_{0},yA_{0})$.

Now, since $xs<x$, we have
$$xsyA_{0}\subset \{\mu\in V\mid n_{s}<\mathcal{h}\mu,\check{\alpha}_{s}\mathcal{i}\}.$$
Therefore, for all $m\leq n_{s}$, we have $H_{\alpha_{s},m}\notin H(yA_{0},xsyA_{0})$. Thus, there is no hyperplane parallel to $H_{s}$ in $H_{xs,y}$, as required.
\end{proof}

We keep the setting of the previous lemma. We denote by $\si_{s}$ the reflection with fixed point set $H_{s}$. Assume that $I_{xs,y}\neq \emptyset$ and let $i\in I_{xs,y}$. Recall that $U_{i}(yA_{0})$ is the unique strip of direction $i$ which contains $yA_{0}$. Since $i\in I_{xs,y}$ we have
$$A_{0}\not\subset U_{i}(yA_{0})\text{ and } xsyA_{0}\not\subset U_{i}(yA_{0}).$$
One can see that one and only one of the hyperplanes which defines $U_{i}(yA_{0})$ lies in $H_{xs,y}$. We denote by $\Hi$ this hyperplane.

Let $H\in H_{xs,y}$. By the previous lemma we know that $H$ is not parallel to $H_{s}$. Consider the 4 connected components of $V-\{H,H_{s}\}$. We denote by $E_{A_{0}}$, $E_{yA_{0}}$, $E_{syA_{0}}$ and $E_{xsyA_{0}}$ the connected component which contains, respectively, $A_{0}$, $yA_{0}$, $syA_{0}$ and $xsyA_{0}$. Assume that $(H)\si_{s}\neq H$. Then, we have either
$$(H)\si_{s}\cap E_{yA_{0}}\neq\emptyset \text{ and }(H)\si_{s}\cap E_{A_{0}}\neq\emptyset$$
or
$$(H)\si_{s}\cap E_{xsyA_{0}}\neq\emptyset \text{ and } (H)\si_{s}\cap E_{syA_{0}}\neq\emptyset.$$
Furthermore, in the first case, $(H)\si_{s}$ separates $E_{xsyA_{0}}$ and $E_{syA_{0}}$, and, in the second case, $(H)\si_{s}$ separates $E_{yA_{0}}$ and $E_{A_{0}}$. In particular, we have
\begin{align*}
(H)\si_{s}\cap E_{yA_{0}}\neq\emptyset \quad&\Longrightarrow\quad (H)\si_{s}\in H(syA_{0},xsyA_{0})\\
(H)\si_{s}\cap E_{xsyA_{0}}\neq\emptyset\quad &\Longrightarrow\quad (H)\si_{s}\in H(A_{0},yA_{0}).
\end{align*}
Moreover, we see that
\begin{align*}
\label{one}
(H)\si_{s}\cap E_{yA_{0}}\neq\emptyset\quad &\Longrightarrow\quad (H)\si_{s}\in H(syA_{0},xsyA_{0}) \notag\\
&\Longrightarrow \quad (H)\si_{s}\in H(yA_{0},xsyA_{0}),
\end{align*}
since $H_{s}$ is the only hyperplane in $H(yA_{0},syA_{0})$ and $(H)\sigma_{s}\neq H_{s}$.

We will say that $H\in H_{xs,y}$ is of $s$-type $1$ if $(H)\si_{s}\cap E_{yA_{0}}\neq \emptyset$ and of $s$-type 2 if $(H)\si_{s}\cap E_{xsyA_{0}}\neq \emptyset$. 
To sum up, we have
\begin{enumerate}
\item[-] if $H$ is of $s$-type 1 then $(H)\si_{s}\in H(yA_{0},xsyA_{0})$;
\item[-] if $H$ is of $s$-type 2 then $(H)\si_{s}\in H(A_{0},yA_{0})$.
\end{enumerate}
We illustrate this result in Figure 1. Note that if $H, H'\in H_{xs,y}$ are parallel, then they have the same type.

\begin{center}
\begin{pspicture}(-5,3.3)(5,-3.3)
\psset{unit=0.75cm}
\psline(-5.5,0)(-1,0)
\psline(6,0)(1,0)
\psline(-2,3)(-4,-3)
\psline[linestyle=dashed](-2,-3)(-4,3)
\psline(-1.5,1)(-2,1)
\psline(-1.5,1)(-1.75,1.5)
\psline(-2,1)(-1.75,1.5)
\rput(-1.625,0.8){$xsyA_{0}$}
\psline(-1.5,-1)(-2,-1)
\psline(-1.5,-1)(-1.75,-1.5)
\psline(-2,-1)(-1.75,-1.5)
\rput(-1.625,-0.7){$xyA_{0}$}

\psline(-1,-2.5)(-1,-3)
\psline(-1,-2.5)(-0.5,-2.75)
\psline(-1,-3)(-0.5,-2.75)
\rput(-0.7,-3.15){$A_{0}$}

\psline(-5,0.5)(-5,-0.5)
\psline(-4.5,0)(-5,0.5)
\psline(-5,-0.5)(-4.5,0)
\rput(-4.3,0.48){$yA_{0}$}
\rput(-4.3,-0.45){$syA_{0}$}

\rput(-0.7,0){$H_{s}$}
\rput(-1.7,3){$H$}
\rput(-4.8,3){$(H)\si_{s}$}

\rput(-0.7,2){$E_{xsyA_{0}}$}
\rput(-0.7,-2){$E_{A_{0}}$}
\rput(-5,-2){$E_{syA_{0}}$}
\rput(-5,2){$E_{yA_{0}}$}

\rput(-2.9,4){$s$-type 1}

\rput(3.6,4){$s$-type 2}
\psline(4.75,1)(4.25,1)
\psline(4.75,1)(4.5,1.5)
\psline(4.25,1)(4.5,1.5)
\rput(4.625,0.8){$xsyA_{0}$}
\psline(4.75,-1)(4.25,-1)
\psline(4.75,-1)(4.5,-1.5)
\psline(4.25,-1)(4.5,-1.5)
\rput(4.625,-0.72){$xyA_{0}$}

\psline(5.25,-2.5)(5.25,-3)
\psline(5.25,-2.5)(5.75,-2.75)
\psline(5.25,-3)(5.75,-2.75)
\rput(5.55,-3.15){$A_{0}$}

\psline(1.55,0.5)(1.55,-0.5)
\psline(2.05,0)(1.55,0.5)
\psline(1.55,-0.5)(2.05,0)
\rput(2.2,0.5){$yA_{0}$}
\rput(2.2,-0.45){$syA_{0}$}

\rput(6.25,0){$H_{s}$}
\rput(5.3,3){$(H)\si_{s}$}
\rput(2.25,3){$H$}

\psline[linestyle=dashed](4.55,3)(2.55,-3)
\psline(4.55,-3)(2.55,3)

\rput(6.25,2){$E_{xsyA_{0}}$}
\rput(6.25,-2){$E_{A_{0}}$}
\rput(1.55,-2){$E_{syA_{0}}$}
\rput(1.55,2){$E_{yA_{0}}$}

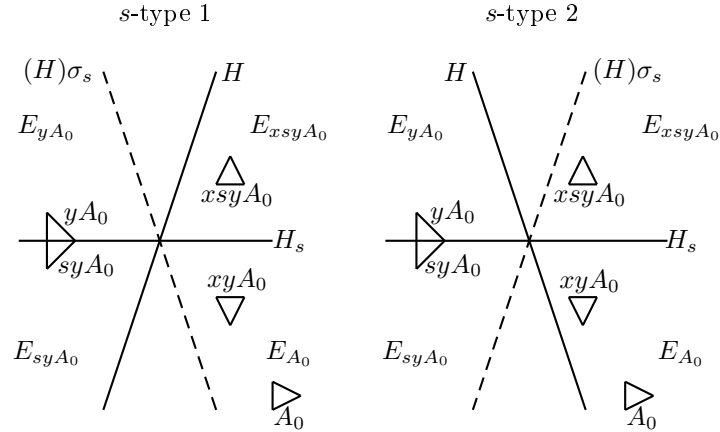
\rput(0.5,-4){FIGURE 1. $s$-type 1 and $s$-type 2 hyperplanes}

\end{pspicture}
\end{center}

\begin{Lem}
\label{Lem3} 
Let $x,y\in W$ and $s\in S$ be such that $xs<x$ and $sy<y$. Let $H_{s}$ be the unique hyperplane which separates $yA_{0}$ and $syA_{0}$ and let $\si_{s}$ be the corresponding reflection. The following holds.
\begin{enumerate}
\item[a)] Let $H\in\mathcal{F}$. We have
\begin{eqnarray*}
H\in H(yA_{0},xsyA_{0})\Rightarrow (H)\si_{s}\in H(yA_{0},xyA_{0}).
\end{eqnarray*}
\item[b)] Let $H\in H_{xs,y}$ be of $s$-type 1; then $H\in H_{x,y}$.
\item[c)] Let $H\in H_{xs,y}$ be of $s$-type 2; then $(H)\si_{s}\in H_{x,y}$.
\item[d)] Let $H\in H_{xs,y}$ such that $(H)\si_{s}=H$; then $H\in H_{x,y}$.
\end{enumerate}
\end{Lem}
\begin{proof}
We prove (a). Let $H\in H(yA_{0},xsyA_{0})$. Then $(H)\si_{s}$ separates $yA_{0}\sigma_{s}$ and $xsyA_{0}\sigma_{s}$. But we have
$$yA_{0}\sigma_{s}=syA_{0}\quad\text{and}\quad xsyA_{0}\sigma_{s}=xssyA_{0}=xyA_{0}.$$
Since $H\neq H_{s}$, we have $(H)\si_{s}\neq H_{s}$ and this implies that $(H)\si_{s}$ separates $yA_{0}$ and $xyA_{0}$.

We prove (b).  We have $H\in H_{xs,y}=H(A_{0},yA_{0})\cap H(yA_{0},xsyA_{0})$. The hyperplane $H$ is of $s$-type 1 thus $(H)\si_{s}\in H(yA_{0},xsyA_{0})$. Using (a) we see that $H\in H(yA_{0},xyA_{0})$. Therefore, $H\in H_{x,y}$.

We prove (c). Since $H$ is of $s$-type 2 we have $(H)\si_{s}\in H(A_{0},yA_{0})$. Moreover, $H\in H(yA_{0},xsyA_{0})$ thus, using (a), we see that $(H)\si_{s}\in H(yA_{0},xyA_{0})$. Therefore, $(H)\si_{s}\in H_{x,y}$.

We prove (d). Using (a), we see that $(H)\si_{s}=H\in H(yA_{0},xyA_{0})$ and since $H\in H_{xs,y}\subset H(A_{0},yA_{0})$, we get $H\in H_{x,y}$.\\
\end{proof}

\begin{Lem}
\label{Lem4}
Let $x,y\in W$ and $s\in S$ be such that $xs<x$ and $sy<y$. Let $H_{s}$ be the unique hyperplane which separates $yA_{0}$ and $syA_{0}$. There is an injective map $\varphi$ from $I_{xs,y}$ to $I_{x,y}-\{ \overline{H_{s}} \}$.
\end{Lem}
\begin{proof}
Let $\si_{s}$ be the reflection with fixed point set $H_{s}$. If $I_{xs,y}=\emptyset$ then the result is clear. We assume that $I_{xs,y}\neq\emptyset$. We define $\varphi$ as follows. 
\begin{enumerate}
\item If $(\Hi)\si_{s}\in H(A_{0},yA_{0})$ then set $\varphi(i)=\si_{s}(i)$;
\item set $\varphi(i)=i$ otherwise.
\end{enumerate}
We need to show that $\varphi(i)\in I_{x,y}-\{\overline{H_{s}}\}$. The fact that $\varphi(i)\neq \overline{H_{s}}$ is a consequence of Lemma \ref{Lem2}, where we have seen that $\overline{H_{s}}\notin I_{xs,y}$. Indeed, since $\varphi(i)$ is either $i$ or $\sigma_{s}(i)$ and $i\neq \overline{H_{s}}$ we cannot have $\varphi(i)=\overline{H_{s}}$. \\
Let $i\in I_{xs,y}$ be such that $\si_{s}(\Hi)\in H(A_{0},yA_{0})$. By Lemma \ref{Lem3} (a), we have $\si_{s}(\Hi)\in H(yA_{0},xyA_{0})$. It follows that $\si_{s}(\Hi)\in H_{x,y}$ and $\si_{s}(i)\in I_{x,y}$ as required.\\
Let $i\in I_{xs,y}$ be such that $\si_{s}(\Hi)\notin H(A_{0},yA_{0})$. Then $\Hi$ is of $s$-type 1. By the previous lemma we have $\Hi\in H_{x,y}$ and $i\in I_{x,y}$. \\
We show that $\varphi$ is injective. Let $i\in I_{xs,y}$ be such that $\varphi(i)=\si_{s}(i)$ and assume that $\si_{s}(i)\in I_{xs,y}$. We have 
$$(U_{i}(yA_{0}))\si_{s}=U_{\si_{s}(i)}(syA_{0})=U_{\si_{s}(i)}(yA_{0})$$
and $(\Hi)\si_{s}$ is one of the hyperplane which defines $U_{\si_{s}(i)}(yA_{0})$. Furthermore since $(\Hi)\si_{s}\in H(A_{0},yA_{0})$ we must have $(\Hi)\si_{s}=H^{(\si_{s}(i))}$. It follows that $(H^{(\si_{s}(i))})\si_{s}\in H(A_{0},yA_{0})$ and $\varphi(\si_{s}(i))=i$. The result follows.
\end{proof}

\begin{Lem}
\label{Lem5}
Let $x,y\in W$ and $s\in S$ such that $xs<x$ and $sy<y$. Let $H_{s}$ be the unique hyperplane which separates $yA_{0}$ and $syA_{0}$. We have
$$c_{xs,y}\leq c_{x,y}-c_{x,y}(\overline{H_{s}}).$$
\end{Lem}
\begin{proof}
Let $\varphi$ be as in the proof of the previous lemma. We keep the same notation. If $I_{xs,y}=\emptyset$ then the result is clear, thus we may assume that $I_{xs,y}\neq 0$.

First assume that $W$ is not of type $\tilde{C}_{r}$ ($r\geq 2$) or $\tilde{A}_{1}$. Then any two parallel hyperplanes have the same weight, therefore we obtain, for $i\in I_{xs,y}$
$$c_{xs,y}(i)=c_{\Hi}.$$ 
Moreover, since $c_{H}=c_{(H)\si}$ for all $H\in\mathcal{F}$ and $\sigma\in\Omega$, one can see that
$$c_{xs,y}(i)=c_{x,y}(\varphi(i)),$$
and the result follows using Lemma \ref{Lem4}.

Now, assume that $W$ is of type $\tilde{C}_{r}$, with graph and weight function given by
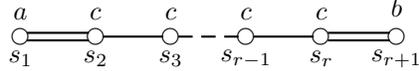
\begin{center}
\begin{pspicture}(0,0)(7,1)
\rput(1,0.5){\circle{0.2}}
\psline(0.99,0.55)(1.81,0.55)
\psline(0.99,0.45)(1.81,0.45)
\rput(2,0.5){\circle{0.2}}
\psline(2,0.5)(2.8,0.5)
\rput(3,0.5){\circle{0.2}}
\psline[linestyle=dashed](3,0.5)(3.8,0.5)
\rput(4,0.5){\circle{0.2}}
\psline(4,0.5)(4.8,0.5)
\rput(5,0.5){\circle{0.2}}
\psline(4.99,0.55)(5.81,0.55)
\psline(4.99,0.45)(5.81,0.45)
\rput(6,0.5){\circle{0.2}}
\rput(0.9,0.2){$s_{1}$}
\rput(0.9,0.8){$a$}
\rput(1.9,0.2){$s_{2}$}
\rput(1.9,0.8){$c$}
\rput(2.9,0.2){$s_{3}$}
\rput(2.9,0.8){$c$}
\rput(3.9,0.2){$s_{r-1}$}
\rput(3.9,0.8){$c$}
\rput(4.9,0.2){$s_{r}$}
\rput(4.9,0.8){$c$}
\rput(5.9,0.2){$s_{r+1}$}
\rput(5.9,0.88){$b$}
\end{pspicture}
\end{center}
In \cite{Bremke}, Bremke proved that the only case where two parallel hyperplanes $H$, $H'$ do not have the same weight is when one of them, say $H$, supports a face of type $s_{1}$ and $H'$ supports a face of type $s_{r+1}$. 

If $a=b$, then parallel hyperplanes have the same weight and we can conclude as before.

Now assume that $a>b$. Let $i\in\overline{\mathcal{F}}$ be such that not all the hyperplanes with direction $i$ have the same weight. Let $H=H_{\alpha,n}$ be a hyperplane with direction $i$ and weight $a$. Then $H_{\alpha,n-1}$ and $H_{\alpha,n+1}$ have weight $b$ because otherwise all the hyperplanes with direction $i$ would have weight $a$.

Let $i\in I_{xs,y}$. If no hyperplane of direction $i$ supports a face of type $s_{1}$ or $s_{r+1}$ then, as before, we can conclude that
$$c_{xs,y}(i)=c_{x,y}(\varphi(i)).$$

Now, in order to prove that $c_{xs,y}\leq c_{x,y}-c_{x,y}(\overline{H_{s}})$, the only problem which may appear is when there exists $i\in I_{xs,y}$ such that
$$c_{xs,y}(i)=a\quad\text{and}\quad c_{x,y}(\varphi(i))=b.$$
Fix such a $i\in I_{xs,y}$. We claim that
\begin{enumerate}
\item $H^{(i)}$ is of $s$-type 1 and $(H^{(i)})\si_{s}\in H(A_{0},yA_{0})$;
\item $\si_{s}(i)\in I_{xs,y}$, $\varphi(\si_{s}(i))=i$ and
$$c_{xs,y}(\si_{s}(i))=b\quad\text{and}\quad c_{x,y}(i)=a.$$
\end{enumerate}
We prove (1). Let $j\in I_{xs,y}$ be such that $\Hj$ is of $s$-type 2. Then $(\Hj)\si_{s}\in H(A_{0},yA_{0})$ and $\varphi(j)=\si_{s}(j)$. Let $H\in H_{xs,y}$ be such that $\overline{H}=j$. Then $H$ is also of $s$-type 2 and $(H)\si_{s}\in H_{x,y}$ (see Lemma \ref{Lem3}(c)). It follows that $c_{x,y}(\va(j))\geq c_{xs,y}(j)$.\\
Let $j\in I_{xs,y}$ be such that $(\Hj)\si_{s}=\Hj$. Then $(\Hj)\si_{s}\in H(A_{0},yA_{0})$ and $\varphi(j)=j$. Let $H\in H_{xs,y}$ be such that $\overline{H}=j$. Then $(H)\si_{s}=H$ and $H\in H_{x,y}$ (see Lemma \ref{Lem3}(d)). It follows that $c_{x,y}(\va(j))\geq c_{xs,y}(j)$.\\
 Finally let $j\in I_{xs,y}$ be such that $\Hj$ is of $s$-type 1 and $\Hj\notin H(A_{0},yA_{0})$. Then $\varphi(j)=j$. Let $H\in H_{xs,y}$ be such that $\overline{H}=j$. Then $H$ is also of $s$-type 1 and $H\in H_{x,y}$ (see Lemma \ref{Lem3}(b)). It follows that $c_{x,y}(\va(j))\geq c_{xs,y}(j)$.\\
Thus since $c_{x,y}(\va(i))< c_{xs,y}(i)$, we get that $\Hi$ is of $s$-type 1 and $(H^{(i)})\si_{s}\in H(A_{0},yA_{0})$. \\

We prove (2). We know that $\Hi$ is of $s$-type 1 and $(\Hi)\si_{s}\in H(A_{0},yA_{0})$. Thus $(\Hi)\si_{s}\in H_{x,y}$ and $\va(i)=\si_{s}(i)$. In particular, since $c_{x,y}(\varphi(i))=b$, we must have $c_{(\Hi)\si_{s}}=b$, which implies that $c_{\Hi}=b$. \\
Since $\Hi$ is of $s$-type 1 we have $(\Hi)\si_{s}\in H(yA_{0},xsyA_{0})$ which implies that $(\Hi)\si_{s}\in H_{xs,y}$. Thus $\si_{s}(i)\in I_{xs,y}$. Arguing as in the proof of Lemma \ref{Lem4}, we obtain $(\Hi)\si_{s}=H^{(\si_{s}(i))}$ and $\va(\si_{s}(i))=i$.\\
Let $\al\in\Phi^{+}$ and $n\in\nZ$ be such that $\Hi=H_{\al,n}$. Since $c_{xs,y}(i)=a$, one can see that one of the hyperplanes $H_{\al,n-1}$, $H_{\al,n+1}$ lies in $H_{xs,y}$. We denote this hyperplane by $H$. Note that $c_{H}=a$ thus, since $c_{x,y}(\si_{s}(i))=b$, we cannot have $(H)\si_{s}\in H_{x,y}$. Both hyperplanes $(H)\si_{s}$ and $(\Hi)\si_{s}$ separate $yA_{0}$ and $xyA_{0}$ but only $(\Hi)\si_{s}$ lies in $H_{x,y}$. This implies that $A_{0}$ lies in the strip defined by $(H)\si_{s}$ and $(\Hi)\si_{s}$. Since $(\Hi)\si_{s}=H^{(\si_{s}(i))}$ this shows that the only hyperplane of direction $\si_{s}(i)$ which lies in $H_{xs,y}$ is $H^{(\si_{s}(i))}$. Thus we have $c_{xs,y}(\si_{s}(i))=b$. Moreover $\va(\si_{s}(i))=i$ and $H$ is of $s$-type 1, thus $H\in H_{x,y}$ (see Lemma \ref{Lem3} (b)) and $c_{x,y}(i)=a$, as required.\\

Let $I_{>}$ be the subset of $I_{xs,y}$ which consists of the directions $i$ such that $c_{xs,y}(i)=a$ and $c_{x,y}(\varphi(i))=b$. Using (1) and (2), we see that the set $\si_{s}(I_{>})$ is a subset of $I_{xs,y}$ such that for all $i\in  \si_{s}(I_{>})$ we have $c_{xs,y}(i)=b$ and $c_{x,y}(\varphi(i))=a$. Therefore we can conclude that $c_{xs,y}\leq c_{x,y}-c_{x,y}(\ov{H_{s}})$ in the case where $W$ is of type $\tC_{r}$ ($r\geq 2$).\\

In the case where $W$ is of type $\tA_{1}$, the result is clear, since we always have $I_{xs,y}=\emptyset$. The lemma is proved.
\end{proof}

\begin{proof}[Proof of Theorem \ref{bound}] 
Let $x,y\in W$ and
$$T_{x}T_{y}=\underset{z\in W}{\sum}f_{x,y,z}T_{z}\ \ \text{ where $f_{x,y,z}\in\cA$}.$$
We want to prove that the degree of $f_{x,y,z}$ in $v$ is less than or equal to $c_{x,y}$. We proceed by induction $\ell(x)+\ell(y)$. \\
If $\ell(x)+\ell(y)=0$ the result is clear. \\
If $c_{x,y}=0$ then $H_{x,y}=\emptyset$ and $xy=x.y$. Thus $T_{x}T_{y}=T_{xy}$ and the result follows.\\
We may assume that $H_{x,y}\neq \emptyset$, which implies that $\ell(x)>0$ and $\ell(y)>0$. 
Let $x=s_{k}\ldots s_{1}$ be a reduced expression of $x$. There exists $1\leq i\leq k$ such that 
$$\ell(s_{i-1}\ldots s_{1}y)=\ell(y)+i-1 \text{ and }  s_{i}s_{i-1}\ldots s_{1}y<s_{i-1}\ldots s_{1}y.$$ 
Let $x_{0}=s_{k}\ldots s_{i}$ and $y_{0}=s_{i-1}\ldots s_{1}y$. Let $H_{s_{i}}$ be the unique hyperplane which separates $y_{0}A_{0}$ and $s_{i}y_{0}A_{0}$. Note that $c_{H_{s_{i}}}=L(s_{i})$. We have
\begin{align*}
T_{x}T_{y}&=T_{x_{0}}T_{y_{0}}
\end{align*}
Using Lemma \ref{Lem1}, we obtain $c_{x,y}=c_{x_{0},y_{0}}$. We have
\begin{align*}
T_{x_{0},y_{0}}&=T_{s_{k}\ldots s_{i+1}}T_{s_{i}}T_{y_{0}}\\
&=T_{s_{k}\ldots s_{i+1}}(T_{s_{i}y_{0}}+\xi_{s_{i}}T_{y_{0}})\\
&=T_{s_{k}\ldots s_{i+1}}T_{s_{i}y_{0}}+\xi_{s_{i}}T_{s_{k}\ldots s_{i+1}}T_{y_{0}}\\
&=T_{x_{0}s_{i}}T_{s_{i}y_{0}}+\xi_{s_{i}}T_{x_{0}s_{i}}T_{y_{0}}
\end{align*}
By induction, $T_{x_{0}s_{i}}T_{s_{i}y_{0}}$ is an $\mathcal{A}$-linear combination of $T_{z}$ with coefficients of degree less than or equal to $c_{x_{0}s_{i},s_{i}y_{0}}$. Using Lemma \ref{lem0}, we have $c_{x_{0}s_{i},s_{i}y_{0}}\leq c_{x_{0},y_{0}}=c_{x,y}$.\\
By induction, $T_{x_{0}s_{i}}T_{y_{0}}$ is an $\mathcal{A}$-linear combination of $T_{z}$ with coefficients of degree less than or equal to $c_{x_{0}s_{i},y_{0}}$. Therefore the degree of the polynomials occurring in $\xi_{s_{i}}T_{x_{0}s_{i}}T_{y_{0}}$ is less than or equal to $L(s_{i})+c_{x_{0}s_{i},y_{0}}$. Applying Lemma \ref{Lem5} to $x_{0}$ and $y_{0}$ we obtain
$$c_{x_{0}s_{i},y_{0}}\leq c_{x_{0},y_{0}}-c_{x_{0},y_{0}}(\ov{H_{s_{i}}})$$

Since $c_{x_{0},y_{0}}(\overline{H_{s_{i}}})\geq c_{H_{s_{i}}}=L(s_{i})$ we obtain
$$L(s_{i})+c_{x_{0}s_{i},y_{0}}\leq c_{x_{0},y_{0}} =c_{x,y}.$$
The theorem is proved.

\end{proof}

\section{The lowest two-sided cell}

\begin{what}{\bf Kazhdan-Lusztig cells.}
Let $(W,S)$ be a Coxeter group and $L$ a weight function on $W$. Let $\cA=\mathbb{Z}[v,v^{-1}]$ and $\cH$ be the generic Iwahori-Hecke algebra corresponding to $(W,S)$ with parameters $\{L(s)|s\in S\}$.

Let $a\mapsto \overline{a}$ be the involution of $\cA$ which takes $v^{n}$ to $v^{-n}$ for all $n\in\nZ$.
We can extend it to a ring involution from $\cH$ to itself by the formula
$$\overline{\underset{w\in W}{\sum}a_{w}T_{w}}=\underset{w\in W}{\sum}\overline{a}_{w}T^{-1}_{w^{-1}}\ , \text{ where $a_{w}\in \cA$}.$$
Let $\cA_{\leq 0}=\nZ[v^{-1}]$ and $\cA_{<0}=v^{-1}\nZ[v^{-1}]$. For $w\in W$ there exists a unique element $C_{w}\in\cH$ such that
$$\overline{C}_{w}=C_{w} \text{ and } C_{w}=T_{w}+\underset{y<w}{\underset{y\in W}{\sum}}P_{y,w}T_{w}, $$
where $P_{y,w}\in \cA_{<0}$ for $y<w$. In fact, the set $\{C_{w},w\in W\}$ forms a basis of $\cH$, known as the Kazhdan-Lusztig basis. The elements $P_{y,w}$ are called the Kazhdan-Lusztig polynomials. We set $P_{w,w}=1$ for any $w\in W$. 

The Kazhdan-Lusztig left preorder $\leq_{L}$ on $W$ is the relation generated by
$$
\begin{cases}
y\leq_{L} w \text{ if there exists some $s\in S$ such that $C_{y}$ appears with} \\
\text{ a non-zero coefficient in $T_{s}C_{w}$, expressed in the $C_{w}$-basis}
\end{cases}
$$
One can see that
$$\cH C_{w}\subseteq \underset{y\leq_{L} w}{\sum}\cA C_{y}\text{ for any $w\in W$.}$$
The equivalence relation associated to $\leq_{L}$ will be denoted by $\sim_{L}$ and the corresponding equivalence classes are called the left cells of $W$. Similarly, we define $\leq_{R}$, $\sim_{R}$ and  right cells. We say that $x\leq_{LR} y$ if there exists a sequence
$$x=x_{0}, x_{1},..., x_{n}=y$$
such that for all $1\leq i\leq n$ we have $x_{i-1}\leq_{L} x_{i}$ or $x_{i-1}\leq_{R} x_{i}$. We write $\sim_{LR}$ for the associated equivalence relation and the equivalence classes are called two-sided cells. The preorder $\leq_{LR}$ induces a partial order on the two-sided cells of $W$.
\end{what}

\begin{what}{\bf The lowest two-sided cell.}
Let $(W,S)$ be an irreducible affine Weyl group.
In this section, we look at the set
$$c_{0}=\{w\in W\mid w=z'.w_{\lambda}.z,\ z,z'\in W,\ \lambda\in T\}$$ 
where $T$ is the set of special points (see Section 2).
We show that $c_{0}$ is the lowest two-sided cell and we determine the decomposition of $c_{0}$ into left cells.

Recall that, for $\lambda$ a special point, $W_{\lambda}$ is the stabilizer in $W$ of the set of alcoves containing $\lambda$ in their closure, $w_{\lambda}$ is the longest element of $W_{\lambda}$ and $S_{\lambda}=S\cap W_{\lambda}$. In particular, we have
$$sw_{\lambda}<w_{\lambda} \text{ for any $s\in S_{\lambda}$}.$$
For $\lambda\in T$ and $z\in W$ such that $w_{\la}z=w_{\la}.z$, we set
$$N_{\lambda,z}=\{w\in W\mid w=z'.w_{\lambda}.z,\ z'\in W\}.$$
In \cite[Proposition 5.1]{Bremke}, it is shown that $N_{\la,z}$ is included in a left cell.\\
For $\la\in T$, we set
$$M_{\lambda}=\{z\in W\mid w_{\la}z=w_{\lambda}.z,\ sw_{\lambda}z\notin c_{0}\text{ for all $s\in S_{\lambda}$}\}.$$
Following \cite{Shi2}, we choose a set of representatives for the $\Omega$-orbits on $T$ and denote it by $R$. Then
$$c_{0}=\underset{\lambda\in R,\ z\in M_{\lambda}}{\bigcup}N_{\lambda,z}\quad \text{(disjoint union)}$$
It is known (\cite{Shi2}) that this is a union over $|W_{0}|$ terms.

We are now ready to state the main result of this paper.
\begin{Th}
\label{main}
Let $\lambda\in T$ and $z\in M_{\la}$. The set $N_{\la,z}$ is a union of left cells. \\
Furthermore, for $y\in W$ and $w\in N_{\la,z}$, we have
$$y\leq_{L} w\ \Longrightarrow \ y\in N_{\la,z}.$$
\end{Th}
The proof of this theorem will be given in the next section. We now discuss a number
of consequences of Theorem \ref{main}.
\begin{Cor}
\label{lc}
Let $\lambda\in T$ and $z\in M_{\la}$. Then, the set $N_{\la,z}$ is a left cell.
\end{Cor}
\begin{proof}
The set $N_{\la,z}$ is a union of left cells which is included in a left cell (see \cite[Proposition 5.1]{Bremke}). Hence it is a left cell.
\end{proof}
The next step is to prove the following.
\begin{Prop}
The set $c_{0}$ is included in a two-sided cell.
\end{Prop}
\begin{proof}
Let $R=\{\la_{1},...,\la_{n}\}$ be a set of representatives for the $\Omega$-orbits on $T$. For example, if $W$ is of type $\tilde{G}_{2}$ we have $n=1$ and if $W$ is of type $\tilde{B_{n}}$ ($n\geq 3$) we have $n=2$. Set
$$c_{\la_{i}}=\{w\in W|\ w=z'.w_{\la_{i}}.z,\ z,z'\in W \}$$
One can see that
$$c_{0}=\bigcup_{i=1}^{i=n}c_{\la_{i}}$$
and for $1\leq i\leq j\leq n$ we have $c_{\la_{i}}\cap c_{\la_{j}}\neq\emptyset$. Therefore to prove the proposition, it is enough to show that each of the set $c_{\la_{i}}$ is included in a two-sided cell.

Fix $1\leq i\leq n$. Let $w,w' \in c_{\la_{i}}$ and $z,z',y,y'\in W$ be such that $w=z'.w_{\la_{i}}.z$ and $w'=y'.w_{\la_{i}}.y$. Using \cite[Proposition 5.1]{Bremke}, together with its version for right cells, we obtain
$$z'w_{\la_{i}}z\ \sim_{L}\ w_{\la_{i}}z\ \sim_{R} \ w_{\la_{i}}y\ \sim_{R}\ y'w_{\la_{i}}y.$$
The result follows.
\end{proof}

Finally, combining the previous results of Shi, Xi and Bremke with Theorem \ref{main}, we now obtain the following description of the lowest two-sided cell in complete generality.
\begin{Th}
\label{ltsc}
Let $W$ be an irreducible affine Weyl group with associated Weyl group $W_{0}$.
Let
$$c_{0}=\{w\in W|\ w=z'.w_{\la}.z,\ z,z'\in W, \la\in T\}$$
where $T$ is the set of special points. We have
\begin{enumerate}
\item $c_{0}$ is a two-sided cell.
\item $c_{0}$ is the lowest two-sided cell for the partial order on the two-sided cell induced by the preorder $\leq_{LR}$.
\item $c_{0}$ contains exactly $|W_{0}|$ left cells. 
\item The decomposition of $c_{0}$ into left cells is as follows
$$c_{0}=\underset{\la\in R,\ z\in M_{\la}}{\bigcup}N_{\la,z}.$$
\end{enumerate}
\end{Th}
\begin{proof}
We have seen that $c_{0}$ is included in a two-sided cell. Let $w\in c_{0}$ and $y\in W$ such that $y\sim_{LR} w$. In particular we have $y\leq_{LR} w$. We may assume that $y\leq_{L} w$ or $y\leq_{R} w$. We know that 
$$c_{0}=\underset{\la\in R,\ z\in M_{\la}}{\bigcup}N_{\la,z}.$$
Thus $w\in N_{\la,z}$ for some $\la\in R$ and $z\in M_{\la}$. If $y\leq_{L} w$ then, using Theorem \ref{main}, we see that $y\in N_{\la,z}$ and thus $y\in c_{0}$. If $y\leq_{R} w$ then using \cite[\S 8.1]{bible}, we have $y^{-1}\leq_{L} w^{-1}$. But $c_{0}$ is stable by taking the inverse, thus, as before, we see that $y^{-1}\in c_{0}$ and $y\in c_{0}$. This implies that $c_{0}$ is a two-sided cell and that it is the lowest one with respect to $\leq_{LR}$.

By \cite{Shi2}, we know that
$$c_{0}=\underset{\la\in R,\ z\in M_{\la}}{\bigcup}N_{\la,z}$$
is a disjoint union over $|W_{0}|$ terms. By Corollary \ref{lc}, the result follows.
\end{proof}
\end{what}

\section{Proof of Theorem \ref{main}}
We keep the setting of the previous section. For $\la\in T$ we denote by $X_{\la}$ the set of minimal left coset representatives of $W_{\la}$ in $W$, that is
$$X_{\la}:=\{w\in W| \ell(ws)>\ell(w)\text{ for all $s\in S_{\la}$}\}.$$
One can easily check that
$$X_{\la}=\{z\in W| zw_{\la}=z.w_{\la}\} \quad\text{and}\quad X_{\la}^{-1}=\{z\in W| w_{\la}z=w_{\la}.z\}.$$
Let $\la\in T$ and $z\in X_{\la}^{-1}$. For $z'\in W$, we have the following equivalence
$$z'w_{\la}z=z'.w_{\la}.z \Longleftrightarrow z'\in X_{\la}.$$
Indeed, if $z'w_{\la}z=z'.w_{\la}.z$ then we must have $z'w_{\la}=z'.w_{\la}$ and $z'\in X_{\la}$. Conversely, if $z'\in X_{\la}$ then since $z\in X_{\la}^{-1}$ we have $z'w_{\la}z=z'.w_{\la}.z$ (see \cite[Lemma 3.2]{Shi2}). Therefore we see that
$$N_{\la,z}=\{w\in W| w=xw_{\la}z, x\in X_{\la}\}.$$

\begin{what}{\bf Preliminaries.}
\begin{Lem}
\label{endp}
Let $\la\in T$, $z\in M_{\la}$, $y\in X_{\la}$ and $v_{1}<w_{\la}z$. Then, $P_{v_{1},w_{\la}z}T_{y}T_{v_{1}}$ is an $\mathcal{A}$-linear combination of $T_{z'}$ ($z'\in W$) with coefficient in $\mathcal{A}_{<0}$.
\end{Lem}
\begin{proof}
We can write uniquely $v_{1}=w.v'$, where $w\in W_{\lambda}$ and $v'\in X_{\la}^{-1}$ (see \cite[Proposition 2.1.1]{gp}). First, assume that $w=w_{\lambda}$. In that case, we have $yv_{1}=y.v_{1}$ and $T_{y}T_{v_{1}}=T_{yv_{1}}$. Since $P_{v_{1},w_{\la}z}\in\mathcal{A}_{<0}$ the result follows.

Next, assume that $w< w_{\lambda}$. Let $w_{v_{1}}\in W$ such that $w_{v_{1}}w=w_{\lambda}$. The Kazhdan-Lusztig polynomials satisfy the following relation (see \cite[Theorem 6.6.c]{bible})
$$P_{x,w}=v^{-L(s)}P_{sx,w}, \text{ where $x<sx$ and $sw<w$ }.$$
Therefore, one can see that
\begin{enumerate}
\item $P_{v_{1},w_{\la}z}\in v^{-L(w_{v_{1}})}\mathcal{A}_{<0}$ if $w_{\lambda}.v'<w_{\la}z,$
\item $P_{v_{1},w_{\la}z}= v^{-L(w_{v_{1}})}$ if $w_{\lambda}.v'=w_{\la}z.$
\end{enumerate}
Thus, to prove the lemma, it is sufficient to show that the polynomials occurring in the expression of $T_{y}T_{v_{1}}$ in the standard basis are of degree less than or equal to $L(w_{v_{1}})$ in the first case and $L(w_{v_{1}})-1$ in the second case.

Using Theorem \ref{bound}, we know that the degree of these polynomials is less than or equal to $c_{y,v_{1}}$ (for the definition of $c_{y,v_{1}}$, see Section 2.3). 

Let $w_{v_{1}}=s_{n}...s_{m+1}$ and $w=s_{m}...s_{1}$ be reduced expressions, and let $H_{i}$ be the unique hyperplane which separates $s_{i-1}...s_{1}v'A_{0}$ and $s_{i}...s_{1}v'A_{0}$. Note that $c_{H_{i}}=L(s_{i})$.
Let $\lambda'$ be the unique special point contained in the closure of $v'A_{0}$ and $w_{\lambda}v'A_{0}$ (note that $W_{\la'}=W_{\la}$). By definition of $X_{\la}$, $yv_{1}A_{0}$ lies in the quarter $\mathcal{C}$ with vertex $\lambda'$ which contains $v_{1}A_{0}$.

Let $1\leq i\leq m$. Let $\alpha_{i}$ and $k\in\mathbb{Z}$ such that $H_{i}=H_{\alpha_{i},k}$. Assume that $k>0$ (the case $k\leq 0$ is similar). We have $v_{1}A_{0}\in V_{H_{i}}^{+}$. Now, since $\la'$ lies in the closure of $v_{1}A_{0}$ and $\lambda'\in H_{i}$, one can see that
$$k<\mathcal{h}\mu,\check{\alpha_{i}}\mathcal{i}<k+1 \text{ for all $\mu\in v_{1}A_{0}$}.$$
Moreover, $yv_{1}A_{0}\in \mathcal{C}$ implies that
$$k<\mathcal{h}\mu,\check{\alpha_{i}}\mathcal{i} \text{ for all $\mu\in yv_{1}A_{0}$}.$$
From there, we conclude that all the hyperplanes $H_{\alpha_{i},l}$ with $l\leq k$ do not lie in $H(v_{1}A_{0},yv_{1}A_{0})$ and that all the hyperplanes $H_{\alpha_{i},l}$ with $l>k$ do not lie in $H(A_{0},v_{1}A_{0})$. Thus $\overline{H_{i}}\notin I_{y,v_{1}}$ and we have
$$I_{y,v_{1}}\subset \{\overline{H_{m+1}},...,\overline{H_{n}}\},$$
which implies
$$c_{y,v_{1}}\leq \overset{i=n}{\underset{i=m+1}{\sum}} c_{y,v_{1}}(\overline{H_{i}}).$$
Now, if $W$ is not of type $\tilde{C}_{r}$ or $\tilde{A}_{1}$ then any two parallel hyperplanes have same weight and we have
$$
c_{y,v_{1}}(\overline{H_{i}})=
\begin{cases}
0 & \mbox{if } i\notin I_{y,v_{1}},\\
L(s_{i}) &\mbox{otherwise}.
\end{cases}
$$
Thus
$$c_{y,v_{1}}\leq \overset{i=n}{\underset{i=m+1}{\sum}} L(s_{i})=L(w_{v_{1}}),$$
as required in the first case.

Assume that $W$ is of type $\tilde{C}_{r}$ or $\tilde{A}_{1}$. Then, one can see that, since $\la'$ is a special point, we have for all $1\leq i\leq n$, $c_{H_{i}}=c_{\overline{H_{i}}}=L(s_{i})$ and we can conclude as before.

Assume that we are in case 2.  Let $j\in\overline{\mathcal{F}}$. Recall that in \cite{Bremke}, the author defined the strip of direction $j$ as the connected component of
$$V-\underset{c_{H}=c_{j}}{\underset{H\in\mathcal{F}, \overline{H}=j}{\bigcup}}H, \quad\text{where}\quad c_{j}=\underset{H,\ \overline{H}=j}{\max}c_{H}.$$
To avoid confusion, we will call them maximal strips of direction $j$. Let 
$$\mathcal{U}(A)=\underset{U \text{ maximal strip}, A\subset U}{\bigcup}U.$$
Now, $v_{1}=w.v'<w_{\lambda}v'=w_{\la}z$ with $z\in M_{\lambda}$. Recall that
$$M_{\lambda}=\{z\in W\mid w_{\lambda}z=w_{\la}.z,\ sw_{\lambda}z\notin c_{0}\text{ for all $s\in S_{\lambda}$}\}$$
thus $v_{1}=w.v'\notin c_{0}$.
In \cite{Bed}, B\'edard showed (in the equal parameter case) that the lowest two-sided cell $c_{0}$ can be described as follows
$$c_{0}=\{w\in W\mid wA_{0}\not\subset \mathcal{U}(A_{0})\}.$$
In \cite{Bremke}, Bremke proved that this presentation remains valid in the unequal parameter case. Therefore since $v_{1}A_{0}\notin c_{0}$, we have $v_{1}A_{0}\in\mathcal{U}(A_{0})$ and there exists a maximal strip $U$ which contains $A_{0}$ and $v_{1}A_{0}$. Let $1\leq k\leq n$ be such that $U$ is of direction $\ov{H_{k}}$. For $1\leq i\leq m$, the hyperplane $H_{i}$ separates $A_{0}$ and $v_{1}A_{0}$ and $c_{H_{i}}=c_{\ov{H_{i}}}$, thus we must have $k>m$. 

If $W$ is not of type $\tilde{C}_{r}$ or $\tilde{A}_{1}$, then our strips and the strips as defined in \cite{Bremke} are the same. Therefore, since  $A_{0}$ and $v_{1}A_{0}$ lie in $U$, we have $\overline{H_{k}}\notin I_{y,v_{1}}$ and
$$c_{y,v_{1}}\leq \underset{i\neq k}{\overset{i=n}{\underset{i=m+1}{\sum}}} c_{y,v_{1}}(\overline{H_{i}})\leq \underset{i\neq k}{\overset{i=n}{\underset{i=m+1}{\sum}}} L(s_{i})<L(w_{v_{1}}),$$
as required.

Assume that $W$ is of type $\tilde{C}_{r}$ or $\tilde{A}_{1}$. 
First, if all the hyperplanes with direction $\overline{H_{k}}$ have same weight, then we have $\overline{H_{k}}\notin I_{y,v_{1}}$ and we can conclude as before.\\
Assume not, then we must have $c_{H_{k}}=c_{\overline{H_{k}}}$ (since $\lambda'\in H_{k}$) and there is no hyperplane of direction $\overline{H_{k}}$ and maximal weight which separates $A_{0}$ and $v_{1}A_{0}$. Therefore
$$c_{y,v_{1}}\leq \overset{i=n}{\underset{i=m+1}{\sum}}
c_{y,v_{1}}(\overline{H_{i}})<\underset{i\neq k}{\overset{i=n}
{\underset{i=m+1}{\sum}}}c_{y,v_{1}}(\overline{H_{i}})+c_{\overline{H_{k}}}\leq \overset{i=n}
{\underset{i=m+1}{\sum}}L(s_{i})=L(w_{v_{1}}),$$
as required.
\end{proof}
\end{what}

\begin{what}{\bf Proof of Theorem \ref{main}.}
In this section we fix $\la\in T$ and $z\in M_{\la}$. We set $v=w_{\la}z$. 
The following argument is inspired by a paper of Geck \cite{Geck}.
\begin{Lem}
The submodule $\mathcal{M}:=\sg T_{x}C_{v}\mid x\in X_{\la}\sd_{\mathcal{A}}\subset \mathcal{H}$ is a left ideal.
\end{Lem}
\begin{proof} Since $\mathcal{H}$ is generated by $T_{s}$ for $s\in S$, it is enough to check that $T_{s}T_{x}C_{v}\in\mathcal{M}$ for $x\in X_{\la}$. According to Deodhar's lemma (see \cite[Lemma 2.1.2]{gp}), there are three cases to consider
\begin{enumerate}
\item $sx\in X_{\la}$ and $\ell(sx)>\ell(x)$. Then $T_{s}T_{x}C_{v}=T_{sx}C_{v}\in \cM$ as required.
\item $sx\in X_{\la}$ and $\ell(sx)<\ell(x)$. Then $T_{s}T_{x}C_{v}=T_{sx}C_{v}+(v_{s}-v_{s}^{-1})T_{x}C_{v}\in \cM$ as required.
\item $t:=x^{-1}sx\in S_{\la}$. Then $\ell(sx)=\ell(x)+1=\ell(xt)$. Now, since $tv<v$, we have (see in \cite[\S 5.5, Theorem 6.6.b]{bible})
$$T_{t}C_{v}=v^{L(t)}C_{v}.$$
Thus, we see that
$$T_{s}T_{x}C_{v}=T_{sx}C_{v}=T_{xt}C_{v}=T_{x}T_{t}C_{v}=v^{L(t)}T_{x}C_{v}$$
which is in $\mathcal{M}$ as required.
\end{enumerate}
\end{proof}
Note that $\{T_{y}C_{v}\mid y\in X_{\la}\}$ is an $\cA$-linearly independent subset of $\cH$. Indeed, for $y\in X_{\la}$ we have
\begin{align*}
T_{y}C_{v}&=T_{y}T_{v}+\underset{u<v}{\sum}P_{u,v}T_{y}T_{u}\\
&=T_{yv}+\text{ an $\cA$-linear combination of $T_{w}$ with $\ell(w)<\ell(yv)$.}
\end{align*}
So, by an easy induction on the length and since the $T_{z}$ form a basis of $\cH$, we get the result. It follows that $\{T_{y}C_{v}\mid y\in X_{\la}\}$ is a basis of $\cM$. Following \cite{Geck}, we shall now exhibit another basis of $\cM$.
\begin{Lem}
Let $y\in X_{\la}$, we can write uniquely
$$T_{y^{-1}}^{-1}C_{v}=\underset{x\in X_{\la}}{\sum}\overline{r}_{x,y}T_{x}C_{v}\ \text{ , }\ \ \overline{r}_{x,y}\in \mathcal{A},$$
where $r_{y,y}=1$ and $r_{x,y}=0$ unless $x\leq y$.
\end{Lem}
\begin{proof}
We have
$$T_{y^{-1}}^{-1}=T_{y}+\underset{z<y}{\sum}\overline{R}_{z,y}T_{z}.$$
where $R_{y,v}\in \cA$ (see \cite[\S 4]{bible}).
Now let $z\in W$ be such that $T_{z}$ occurs in the above expression. We can write $z$ uniquely under the form $x.w$ with $w\in W_{\lambda}$ and $x\in X_{\la}$ (see \cite[Proposition 2.1.1]{gp}). Note that $x\leq z<y$. Since $xw=x.w$, we have $T_{z}=T_{x}T_{w}$ and 
$$T_{y^{-1}}^{-1}C_{v}=T_{y}C_{v}+ \text{an $\mathcal{A}$-linear combination of }T_{x}T_{w}C_{v}, \text{where $x<y$}$$
Since $w\in W_{\la}$, we know that $\ell(ww_{\la})=\ell(w_{\la})-\ell(w)$ and $T_{w}C_{v}=v^{L(w)}C_{v}$. Therefore we see that $r_{y,y}=1$ and $r_{x,y}=0$ unless $x\leq y$.
\end{proof}
\begin{Lem}
\label{delta}
Let $x,y\in X_{\la}$. We have
$$\underset{x\leq z\leq y}{\underset{z\in X_{\la}}{\sum}}\overline{r}_{x,z}r_{z,y}=\delta_{x,y}. $$
\end{Lem}
\begin{proof}
The proof is similar to the one in \cite{Geck}, once we know that $\{T_{y}C_{v}\mid y\in X_{\la}\}$ is an $\cA$-linearly independent subset of $\cH$.
\end{proof}
\begin{Prop}
\label{nb}
For any $y\in X_{\la}$, we have
$$C_{yv}=T_{y}C_{v}+\underset{x\in X_{\la},\ x<y}{\sum}p^{*}_{x,y}T_{x}C_{v}\ \ \text{where} \ p^{*}_{x,y}\in \mathcal{A}_{<0}.$$
\end{Prop}
\begin{proof} 
Fix $y\in X_{\la}$ and consider a linear combination
$$\tilde{C}_{yv}:=\underset{x\in X_{\la}, x\leq y}{\sum}p^{*}_{x,y}T_{x}C_{v}$$ 
where $p^{*}_{y,y}=1$ and $p^{*}_{x,y}\in\mathcal{A}_{<0}$ if $x<y$. 

Our first task is to show that the $p^{*}$-polynomials can be chosen such that $\tilde{C}_{yv}=\overline{\tilde{C}_{yv}}$. In order to do so, we proceed as in the proof of the existence of the Kazhdan-Lusztig basis in \cite{Lus2}. We set up a system of equations with unknown the $p^{*}$-polynomials and then use an inductive argument to show that this system has a unique solution.

As in \cite{Geck}, one sees that the condition $\overline{\tilde{C}}_{yv}=\tilde{C}_{yv}$  is equivalent to
$$\underset{x\in X_{\la},z\leq x\leq y}{\sum}\overline{p}^{*}_{x,y}\overline{r}_{z,x}=p^{*}_{z,y}\ \ \text{for all $z\in X_{\la}$ such that $z\leq y$}.$$
In other words, the coefficient $p^{*}_{x,y}$ must satisfy
\begin{align*}
p^{*}_{y,y}&=1, &(1)\\
\overline{p}^{*}_{x,y}-p^{*}_{x,y}&=\underset{x<z\leq y}{\underset{z\in X_{\la}}{\sum}}r_{x,z}p^{*}_{z,y} \text{ for $x\in X_{\la}$, $x<y$. } &(2)
\end{align*}
We now consider (1) and (2) as a system of equations with unknown $p^{*}_{x,y}$ ($x\in X_{\la}$). We solve it by induction. Let $x\in X_{\la}$ be such that $x\leq y$. If $x=y$, we set $p^{*}_{y,y}=1$, so (1) holds.\\
Now assume that $x<y$. Then, by induction, $p^{*}_{z,y}$ ($z\in X_{\la}$) are known for all $x<z\leq y$ and they satisfy
$p^{*}_{z,y}\in \mathcal{A}_{<0}$ if $z<y$. In other words, the right-hand side of (2) is known, we denote it by $a\in\mathcal{A}$.

Using Lemma \ref{delta} and the same argument as in \cite{Geck} yields $\overline{a}=-a$. Thus the identity $\overline{p}^{*}_{x,y}-p^{*}_{x,y}=a$ together with the condition that $p^{*}_{x,y}\in\mathcal{A}_{<0}$ uniquely determine $p^{*}_{x,y}$.

We have proved that the coefficient $p^{*}_{x,y}$ can be chosen such that $\tilde{C}_{yv}$ is fixed by the involution $h\rightarrow \overline{h}$. Furthermore, we have
\begin{align*}
\tilde{C}_{yv}&=T_{y}C_{v}+\underset{x\in X_{\la}}{\underset{x<y}{\sum}}p^{*}_{x,y}T_{x}C_{v}\\
&= T_{y}(T_{v}+\underset{v_{1}<v}{\sum}P_{v_{1},v}T_{v_{1}})+ \underset{x\in X_{\la}}{\underset{x<y}{\sum}}p^{*}_{x,y}T_{x}(\underset{v_{1}\leq v}{\sum}P_{v_{1},v}T_{v_{1}})\\
&=T_{y}T_{v}+\underset{v_{1}<v}{\sum}P_{v_{1},v}T_{y}T_{v_{1}}+ \underset{x\in X_{\la}}{\underset{x<y}{\sum}}\underset{v_{1}\leq v}{\sum}p^{*}_{x,y}P_{v_{1},v}T_{x}T_{v_{1}}.
\end{align*}
Now, in the above expression, the elements of the form $P_{v_{1},v}T_{x}T_{v_{1}}$, with $x\in X_{\la}$ and $v_{1}<v$, can give some coefficient in $\mathcal{A}_{\geq 0}$ (compare to the situation in \cite{Geck}). However, by Lemma \ref{endp}, we get
$$\tilde{C}_{yv}=T_{yv}+\text{ an $\mathcal{A}_{<0}$-linear combination of $T_{z}$ with $\ell(z)<\ell(yv)$}$$
and then by definition of the Kazhdan-Lusztig basis, we conclude that $\tilde{C}_{yv}=C_{yv}$ as required.
\end{proof}

\begin{Cor}
\label{F}
We have
$$\mathcal{M}=\sg C_{yv}\mid y\in X_{\la}\sd_{\cA}.$$
\end{Cor}
\begin{proof} Let $y\in X_{\la}$. By Proposition \ref{nb} we have
$$C_{yv}=T_{y}C_{v}+\underset{x<y}{\underset{x\in X_{\la}}{\sum}}p^{*}_{x,y}T_{x}C_{v}.$$
Thus $C_{yv}\in\mathcal{M}$. Now, by a straightforward induction, we see that
$$T_{y}C_{v}=C_{yv}+\text{ a $\mathcal{A}$-linear combination of $C_{xv}$},$$
which yields the required assertion.
\end{proof}

We can now prove Theorem \ref{main}.\\
Recall that
$$N_{\lambda,z}=\{w\in W\mid w=xv,\ x\in X_{\la}\}.$$
Let $y\in X_{\la}$. Let $z\in W$ such that $z\leq_{L} yv$. We want to show that $z\in N_{\lambda,z}$. To prove the theorem it is enough to consider
 the case where $C_{z}$ appears with a non-zero coefficient in $T_{s}C_{yv}$ for some $s\in S$.
 By Corollary \ref{F}, $C_{yv}\in\mathcal{M}$. Since $\mathcal{M}$ is a left ideal, we know that 
 $T_{s}C_{yv}\in \mathcal{M}$. Using Corollary \ref{F} once more, we obtain
$$T_{s}C_{yv}=\underset{y_{1}\in X_{\la}}{\sum}a_{y_{1},yv}C_{y_{1}v}$$
and this expression is the expression of $T_{s}C_{yv}$ in the Kazhdan-Lusztig basis of $\mathcal{H}$.  We assumed that $C_{z}$ appears with non zero coefficient in that expression, therefore there exists $y_{1}\in X_{\la}$ such that $z=y_{1}v$, and $z\in N_{\lambda,z}$ as required. The theorem follows.

\end{what}

\end{document}